\documentclass[12pt]{article}
\usepackage{amssymb}

\newtheorem{definition}{Definition}[section]
\newtheorem{theorem}[definition]{Theorem}

\newtheorem{lemma}[definition]{Lemma}

\newtheorem{remark}[definition]{Remark}
\newtheorem{corollary}[definition]{Corollary}

\newtheorem{note}[definition]{Note}
\newtheorem{notation}[definition]{Notation}
\newtheorem{assumption}[definition]{Assumption}
\def\R{\mathbb R}
\def\C{\mathbb C}

\newcommand{\beast}{\begin{eqnarray*}}
\newcommand{\eeast}{\end{eqnarray*}}

\begin{document}

\title{\bf A duality between pairs of split decompositions for a $Q$-polynomial distance-regular graph}

\author{
Joohyung Kim{\footnote{
Department of Mathematics, 
University of Wisconsin,
480 Lincoln Drive,
Madison WI,
53706-1388 USA
}}}

\date{}

\maketitle

\begin{abstract}
Let $\Gamma$ denote a $Q$-polynomial distance-regular graph with diameter $D \geq 3$ and standard module $V$. Recently Ito and Terwilliger introduced four direct sum decompositions of $V$; we call these the $(\mu,\nu)$--{\it split decompositions} of $V$, where $\mu, \nu \in \lbrace \downarrow, \uparrow \rbrace$. In this paper we show that the ($\downarrow,\downarrow$)--split decomposition and the ($\uparrow,\uparrow$)--split decomposition are dual with respect to the standard Hermitian form on $V$. We also show that the ($\downarrow,\uparrow$)--split decomposition and the ($\uparrow,\downarrow$)--split decomposition are dual with respect to the standard Hermitian form on $V$.

\end{abstract}

\section{Introduction}

We consider a distance-regular graph $\Gamma$ with vertex set $X$ and diameter $D \geq 3$ (see Section 4 for formal definitions). We assume that $\Gamma$ is $Q$-polynomial with respect to the ordering $E_0, E_1, \ldots, E_D$ of the primitive idempotents. Let $V$ denote the vector space over $\C$
consisting of column vectors whose coordinates are indexed by $X$ and whose entries are in $\C$.
We call $V$ the {\it standard module}. We endow $V$ with the Hermitian form $\langle \, , \, \rangle$ 
that satisfies
$\langle u,v \rangle = u^t\overline{v}$ for 
$u,v \in V$. We call this form the {\it standard Hermitian form} on $V$. Recently Ito and Terwilliger introduced four direct sum decompositions of $V$ \cite{drgqtet}; we call these the $(\mu,\nu)$--{\it split decompositions} of $V$, where $\mu, \nu \in \lbrace \downarrow, \uparrow \rbrace$. These are defined as follows. Fix a vertex $x \in X$. For $0 \leq i \leq D$ let $E^*_i=E^*_i(x)$ denote the diagonal matrix in $\mbox{Mat}_X(\C)$ that represents the
 projection onto the $i$th subconstituent of
$\Gamma$ with respect to $x$. For $-1\leq i,j\leq D$
we define
\begin{eqnarray*}
V^{\downarrow \downarrow}_{i,j} &=& 
(E^*_0V+\cdots+E^*_iV)\cap (E_0V+\cdots+E_jV),
\\
V^{\uparrow \downarrow}_{i,j} &=& 
(E^*_DV+\cdots+E^*_{D-i}V)\cap (E_0V+\cdots+E_jV),
\\
V^{\downarrow \uparrow}_{i,j} &=& 
(E^*_0V+\cdots+E^*_iV)\cap (E_DV+\cdots+E_{D-j}V),
\\
V^{\uparrow \uparrow}_{i,j} &=& 
(E^*_DV+\cdots+E^*_{D-i}V)\cap (E_DV+\cdots+E_{D-j}V).
\end{eqnarray*}
For $\mu, \nu \in \lbrace \downarrow,
\uparrow \rbrace$ and  $0 \leq i,j\leq D$ we have
$V^{\mu \nu}_{i-1,j} \subseteq V^{\mu \nu}_{i,j}$
and
$V^{\mu \nu}_{i,j-1}  \subseteq V^{\mu \nu}_{i,j}$; therefore
$V^{\mu \nu}_{i-1,j}+
V^{\mu \nu}_{i,j-1} \subseteq V^{\mu \nu}_{i,j}$.
Let ${\tilde V}^{\mu \nu}_{i,j}$
 denote the orthogonal complement of 
$V^{\mu \nu}_{i-1,j}+
V^{\mu \nu}_{i,j-1}$ in $V^{\mu \nu}_{i,j}$ with respect to the standard Hermitian form. By \cite [Lemma 10.3]{drgqtet},\begin{eqnarray*}
V = \sum_{i=0}^D \sum_{j=0}^D {\tilde V}^{\mu \nu}_{i,j}
\qquad \qquad (\mbox{direct sum}).
\end{eqnarray*}
We call the above sum the $(\mu,\nu)$--{\it split decomposition} of $V$ with respect to $x$. We show that with respect to the standard Hermitian form the ($\downarrow,\downarrow$)--split decomposition (resp. ($\downarrow,\uparrow$)--split decomposition) and the ($\uparrow,\uparrow$)--split decomposition (resp. ($\uparrow,\downarrow$)--split decomposition) are dual in the following sense.
 
\begin{theorem}
\label{thm:duality0}
With the above notation, the following (i), (ii) hold for $0 \leq i,j,r,s\leq D$.
\begin{enumerate}
\item ${{\tilde V}^{\downarrow \downarrow}_{i,j}}$ and ${{\tilde V}^{\uparrow \uparrow}_{r,s}}$ are orthogonal unless $i+r=D$ and $j+s=D$.
\item ${{\tilde V}^{\downarrow \uparrow}_{i,j}}$ and ${{\tilde V}^{\uparrow \downarrow}_{r,s}}$ are orthogonal unless $i+r=D$ and $j+s=D$.
\end{enumerate}
\end{theorem}
To prove Theorem \ref{thm:duality0} we use a result about tridiagonal pairs (Theorem \ref{thm:ortho2}) which may be of independent interest. We also use some results about the subconstituent algebra of $\Gamma$.

\section{Tridiagonal pairs}

\noindent We recall the notion of a tridiagonal pair \cite{itt}. We will use the following terms.
Let $V$ denote a vector space over $\mathbb C$ with finite positive dimension.
By a {\it linear transformation} on $V$ we mean a $\mathbb C$-linear map from $V$ to $V$. Let $A$ denote a linear transformation on $V$. By an {\it eigenspace} of $A$ we mean 
a nonzero  subspace of $V$ of the form
\begin{eqnarray*}
\lbrace v \in V \;\vert \;Av = \theta v\rbrace,
\label{eq:defeigspace}
\end{eqnarray*}
where $\theta \in \mathbb C$. We say that $A$ is {\it diagonalizable} whenever
$V$ is spanned by the eigenspaces of $A$. In this case $V$ is the direct sum of the eigenspaces of $A$.

\begin{definition}
\rm
\cite[Definition 1.1]{itt}
\label{def:nonthinlp}
Let $V$ denote
a vector space over $\mathbb C$ with finite positive dimension.
By a {\it tridiagonal pair} (or {\it TD pair}) on $V$ we mean an ordered pair $A,A^*$ of linear transformations on $V$ that satisfy the following four conditions.
\begin{enumerate}
\item $A$ and $A^*$ are both diagonalizable on $V$.
\item There exists an ordering $V_0, V_1,\ldots, V_d$ of the  
eigenspaces of $A$ such that 
\begin{eqnarray*}
A^* V_i \subseteq V_{i-1} + V_i+ V_{i+1} \qquad \qquad (0 \leq i \leq d),
\label{eq:astaractionthreeint}
\end{eqnarray*}
where $V_{-1} = 0$, $V_{d+1}= 0$.
\item There exists an ordering $V^*_0, V^*_1,\ldots, V^*_\delta$ of
the  
eigenspaces of $A^*$ such that 
\begin{eqnarray*}
A V^*_i \subseteq V^*_{i-1} + V^*_i+ V^*_{i+1} \qquad \qquad (0 \leq i \leq \delta),
\label{eq:aactionthreeint}
\end{eqnarray*}
where $V^*_{-1} = 0$, $V^*_{\delta+1}= 0$.
\item There is no subspace $W$ of $V$ such  that  both $AW\subseteq W$,
$A^*W\subseteq W$, other than $W=0$ and $W=V$.
\end{enumerate}
\end{definition}

\begin{note}
\rm
According to a common 
notational convention $A^*$ denotes the 
conjugate-transpose of $A$. We are not using this convention.
In a tridiagonal pair $A,A^*$ the linear transformations 
$A$ and $A^*$ are arbitrary subject to (i)--(iv) above.
\end{note}

\noindent With reference to Definition \ref{def:nonthinlp}, we have $d=\delta$ \cite[Lemma 4.5]{itt}; we call this common value the {\it diameter} of $A,A^*$. See 
\cite{itt,
shape,
tdanduq} for more information on tridiagonal pairs.

\medskip
\noindent With reference to Definition \ref{def:nonthinlp}, by the construction we have the direct sum decompositions $V= \sum_{i=0}^d {V}_{i}$ and $V= \sum_{i=0}^d {V}^*_{i}$. We now recall four more direct sum decompositions of $V$ called the split decompositions.

\begin{lemma}
{\rm
\cite[Lemma 4.2]{tdanduq}}
\label{lem:spl}
With reference to Definition \ref{def:nonthinlp}, for $\mu,\nu \in \{ \downarrow,\uparrow \}$
we have
\begin{eqnarray*}
V &=& \sum_{i=0}^d {U}^{\mu \nu}_{i} \qquad \qquad {\rm(direct\ sum)},
\end{eqnarray*}
where
\begin{eqnarray*}
U^{\downarrow \downarrow}_{i} &=& 
(V^*_0 +\cdots +V^*_i)\cap (V_0 +\cdots +V_{d-i}),
\\
U^{\uparrow \downarrow}_{i} &=& 
(V^*_{d-i} +\cdots +V^*_d)\cap (V_0 +\cdots +V_{d-i}),
\\
U^{\downarrow \uparrow}_{i} &=& 
(V^*_0 +\cdots +V^*_i)\cap (V_i +\cdots +V_d),
\\
U^{\uparrow \uparrow}_{i} &=& 
(V^*_{d-i} +\cdots +V^*_d)\cap (V_i +\cdots +V_d).
\end{eqnarray*}
\end{lemma}

\section{Hermitian forms}

\noindent
In this section we consider a tridiagonal pair for which the underlying vector space supports a certain Hermitian form. We start with the definition of a Hermitian form. Throughout this section $V$ denotes a vector space over $\mathbb C$ with finite positive dimension. For $\alpha \in \mathbb C$ let $\overline \alpha$ denote the complex conjugate of $\alpha$.

\begin{definition}
\label{def:Herm}
\rm
By a {\it Hermitian form} on $V$ we mean a function $(\,,\,): V \times V \rightarrow \mathbb C$ such that for all $u$, $v$, $w$ in $V$ and all $\alpha \in \mathbb C$,
\begin{enumerate}
\item $(u+v,w)=(u,w)+(v,w)$,
\item $(\alpha u,v)=\alpha (u,v)$,
\item $(v,u)=\overline{(u,v)}$.
\end{enumerate}
\end{definition}

\begin{definition}
\label{def:PD}
\rm
Let $(\,,\,)$ denote a Hermitian form on $V$. By Definition \ref{def:Herm}(iii) we have $(v,v) \in \mathbb R$ for $v \in V$. We say that $(\,,\,)$ is {\it positive definite} whenever $(v,v)>0$ for all nonzero $v \in V$.
\end{definition}

\begin{lemma}
\label{lem:Herm}
Let $(\,,\,)$ denote a positive definite Hermitian form on $V$. Suppose that we are given a linear transformation $A: V \to V$ satisfying
\begin{eqnarray}
 (Au,v) &=& (u,Av) \qquad \qquad \qquad u,v \in V.
\label{eq:line11}
\end{eqnarray}
Then all the eigenvalues of $A$ are in $\mathbb R$.
\end{lemma}

\noindent 
{\it Proof:} Let $\lambda$ denote an eigenvalue of $A$. We show that $\lambda \in \mathbb R$. Since $\mathbb C$ is algebraically closed there exists a nonzero $v \in V$ such that $Av=\lambda v$. By (\ref{eq:line11}) $(Av,v) = (v,Av)$. Evaluating this using Definition \ref{def:Herm}(ii),(iii) we have $(\lambda - \overline{\lambda}) (v,v) =0$. But $(v,v) \neq 0$ since $(\,,\,)$ is positive definite so $\lambda = \overline{\lambda}$. Therefore $\lambda \in \mathbb R$. 
\hfill $\Box$ \\

\begin{assumption}
\label{assumption:trisys}
\rm
Let $A,A^*$ denote a tridiagonal pair on $V$ as in Definition \ref{def:nonthinlp}. For $0 \leq i \leq d$ let $\theta_i$ (resp. $\theta^*_i$) denote the eigenvalue of $A$ (resp. $A^*$) associated with $V_i$ (resp. $V^*_i$). 
We remark that $\theta_0, \theta_1, \ldots, \theta_d$ are mutually distinct and $\theta^*_0, \theta^*_1, \ldots, \theta^*_d$ are mutually distinct. We assume that there exists a positive definite Hermitian form $(\,,\,)$ on $V$ satisfying
\begin{eqnarray}
 (Au,v) &=& (u,Av) \qquad \qquad  u,v \in V,
\label{eq:line12}
\\
(A^{*}u,v) &=& (u,A^{*}v) \qquad \qquad  u,v \in V.
\label{eq:line13}
\end{eqnarray}
\end{assumption}

\begin{lemma}
\label{lem:ortho1}
With reference to Assumption \ref{assumption:trisys}, the following (i), (ii) hold.
\begin{enumerate}
\item The eigenspaces $V_0, V_1,\ldots, V_d$ are mutually orthogonal with respect to $(\,,\,)$.
\item The eigenspaces $V^*_0, V^*_1,\ldots, V^*_d$ are mutually orthogonal with respect to $(\,,\,)$.
\end{enumerate}
\end{lemma}

\noindent 
{\it Proof:} (i) For distinct $i,j$ ($0 \leq i,j\leq d$) and for $u \in V_i$, $v \in V_j$ we show that $(u,v)=0$. By (\ref{eq:line12}) $(Au,v)=(u,Av)$. Evaluating this using Definition \ref {def:Herm}(ii),(iii) we find $(\theta_i - \overline{\theta_j})(u,v)=0$. But $\overline{\theta_j}= \theta_j$ by Lemma \ref {lem:Herm} and $\theta_i \neq \theta_j$ so $(u,v)=0$.\\
(ii) Similar to the proof of (i).
\hfill $\Box$ \\

\begin{theorem}
\label{thm:ortho2}
With reference to Lemma \ref{lem:spl} and Assumption \ref{assumption:trisys}, the following (i), (ii) hold for $0\leq i,j\leq d$ such that $i+j\neq d$.
\begin{enumerate}
\item The subspaces $U^{\downarrow \downarrow}_i$ and $U^{\uparrow \uparrow}_j$ are orthogonal with respect to $(\,,\,)$.
\item The subspaces $U^{\downarrow \uparrow}_i$ and $U^{\uparrow \downarrow}_j$ are orthogonal with respect to $(\,,\,)$.
\end{enumerate}
\end{theorem}

\noindent {\it Proof:} (i) We consider two cases: $i+j< d$ and $i+j> d$.
First suppose that $i+j< d$. By Lemma
\ref{lem:spl}, $U^{\downarrow \downarrow}_i \subseteq V^*_0 +\cdots +V^*_i$ and $U^{\uparrow \uparrow}_j \subseteq V^*_{d-j} +\cdots +V^*_d$. Observe that $V^*_0 +\cdots +V^*_i$ is orthogonal to $V^*_{d-j} +\cdots +V^*_d$ by Lemma \ref{lem:ortho1}(ii) and since $i< d-j$. Therefore $U^{\downarrow \downarrow}_i$ is orthogonal to $U^{\uparrow \uparrow}_j$. Next suppose that $i+j> d$. By Lemma \ref{lem:spl},
 $U^{\downarrow \downarrow}_i \subseteq V_0 +\cdots +V_{d-i}$ and $U^{\uparrow \uparrow}_j \subseteq V_j +\cdots +V_d$. Observe that $V_0 +\cdots +V_{d-i}$ is orthogonal to $V_j +\cdots +V_d$ by Lemma \ref{lem:ortho1}(i) and since $d-i< j$. Therefore $U^{\downarrow \downarrow}_i$ is orthogonal to $U^{\uparrow \uparrow}_j$.\\
(ii) Similar to the proof of (i).
\hfill $\Box$

\section{Distance-regular graphs}
\noindent
In this section 
we review some definitions and basic concepts concerning distance-regular
graphs.
For more background information we refer the reader to 
\cite{bannai}, \cite{bcn}, \cite{godsil} and \cite{terwSub1}.

\medskip
\noindent
Let $X$ denote a nonempty  finite  set.
Let
 $\hbox{Mat}_X(\C)$ 
denote the $\C$-algebra
consisting of all matrices whose rows and columns are indexed by $X$
and whose entries are in $\C$. Let
$V=\C^X$ denote the vector space over $\C$
consisting of column vectors whose 
coordinates are indexed by $X$ and whose entries are
in $\C$.
We observe that $\hbox{Mat}_X(\C)$ acts on $V$ by left multiplication.
We call $V$ the {\it standard module}. We endow $V$ with the Hermitian form $\langle \, , \, \rangle$ 
that satisfies
$\langle u,v \rangle = u^t\overline{v}$ for 
$u,v \in V$,
where $t$ denotes transpose. Observe that $\langle \, , \, \rangle$ is positive definite.
We call this form the {\it standard Hermitian form} on $V$. 
Observe that for $B \in \hbox{Mat}_X(\C)$, 
\begin{equation}
\label{SHF}
\langle Bu, v\rangle =\langle u, {\overline{B}}^{t}v \rangle \qquad \qquad  u,v \in V.
\end{equation} 
For all $y \in X,$ let $\hat{y}$ denote the element
of $V$ with a 1 in the $y$ coordinate and 0 in all other coordinates.
Observe that $\{\hat{y}\;|\;y \in X\}$ is an orthonormal basis for $V.$

\medskip
\noindent
Let $\Gamma = (X,R)$ denote a finite, undirected, connected graph,
without loops or multiple edges, with vertex set $X$ and 
edge set
$R$.   
Let $\partial $ denote the
path-length distance function for $\Gamma $,  and set
$D := \mbox{max}\{\partial(x,y) \;|\; x,y \in X\}$.  
We call $D$  the {\it diameter} of $\Gamma $. We say that $\Gamma$ is {\it distance-regular}
whenever for all integers $h,i,j\;(0 \le h,i,j \le D)$ 
and for all vertices $x,y \in X$ with $\partial(x,y)=h,$ the number
\begin{eqnarray*}
p_{ij}^h = |\{z \in X \; |\; \partial(x,z)=i, \partial(z,y)=j \}|
\end{eqnarray*}
is independent of $x$ and $y.$ The $p_{ij}^h$ are called
the {\it intersection numbers} of $\Gamma.$ 

\medskip
\noindent
For the rest of this paper we assume  that $\Gamma$  
is  distance-regular  with diameter $D\geq 3$. 

\medskip
\noindent 
We recall the Bose-Mesner algebra of $\Gamma.$ 
For 
$0 \le i \le D$ let $A_i$ denote the matrix in $\hbox{Mat}_X(\C)$ with
$xy$ entry
$$
{(A_i)_{xy} = \cases{1, & if $\partial(x,y)=i$\cr
0, & if $\partial(x,y) \ne i$\cr}} \qquad (x,y \in X).
$$
We call $A_i$ the $i$th {\it distance matrix} of $\Gamma.$
We abbreviate $A:=A_1$ and call this the {\it adjacency
matrix} of $\Gamma.$ We observe that
(i) $A_0 = I$;
 (ii)
$\sum_{i=0}^D A_i = J$;
(iii)
$\overline{A_i} = A_i \;(0 \le i \le D)$;
(iv) $A_i^t = A_i  \;(0 \le i \le D)$;
(v) $A_iA_j = \sum_{h=0}^D p_{ij}^h A_h \;( 0 \le i,j \le D)
$,
where $I$ (resp. $J$) denotes the identity matrix 
(resp. all 1's matrix) in 
 $\hbox{Mat}_X(\C)$.
 Using these facts  we find
 $A_0,A_1,\ldots,A_D$
is a basis for a commutative subalgebra $M$ of 
$\mbox{Mat}_X(\C)$.
We call $M$ the {\it Bose-Mesner algebra} of $\Gamma$.
It turns out that $A$ generates $M$ \cite[p.~190]{bannai}.
By (\ref{SHF}) and since $A$ is real symmetric,
\begin{equation}
\label{SHF1}
\langle Au, v\rangle =\langle u, Av \rangle \qquad \qquad  u,v \in V.
\end{equation}
By \cite[p.~45]{bcn}, $M$ has a second basis 
$E_0,E_1,\ldots,E_D$ such that
(i) $E_0 = |X|^{-1}J$;
(ii) $\sum_{i=0}^D E_i = I$;
(iii) $\overline{E_i} = E_i \;(0 \le i \le D)$;
(iv) $E_i^t =E_i  \;(0 \le i \le D)$;
(v) $E_iE_j =\delta_{ij}E_i  \;(0 \le i,j \le D)$.
We call $E_0, E_1, \ldots, E_D $  the {\it primitive idempotents}
of $\Gamma$.  

\medskip
\noindent
We recall the eigenvalues of $\Gamma $. Since $E_0,E_1,\ldots,E_D$ form a basis for  
$M$ there exist complex scalars $\theta_0,\theta_1,
\ldots,\theta_D$ such that
$A = \sum_{i=0}^D \theta_iE_i$.
Observe that $AE_i = E_iA =  \theta_iE_i$ for $0 \leq i \leq D$. We call $\theta_i$  the {\it eigenvalue}
of $\Gamma$ associated with $E_i$ $(0 \leq i \leq D)$.
By Lemma \ref{lem:Herm} and 
(\ref{SHF1}) the 
eigenvalues $\theta_0,\theta_1,\ldots,\theta_D$ are
in $\R.$ Observe that
$\theta_0,\theta_1,\ldots,\theta_D$ are mutually distinct 
since $A$ generates $M$.
Observe that
\begin{eqnarray*}
V = E_0V+E_1V+ \cdots +E_DV \qquad \qquad {\rm (orthogonal\ direct\ sum}).
\end{eqnarray*}
For $0 \le i \le D$ the space $E_iV$ is the  eigenspace of $A$ associated 
with $\theta_i$.

\medskip
\noindent 
We now recall the Krein parameters.
Let $\circ $ denote the entrywise product in
$\mbox{Mat}_X(\C)$.
Observe that
$A_i\circ A_j= \delta_{ij}A_i$ for $0 \leq i,j\leq D$,
so
$M$ is closed under
$\circ$. Thus there exist complex scalars
$q^h_{ij}$  $(0 \leq h,i,j\leq D)$ such
that
$$
E_i\circ E_j = |X|^{-1}\sum_{h=0}^D q^h_{ij}E_h
\qquad (0 \leq i,j\leq D).
$$
By \cite[p.~170]{Biggs}, 
$q^h_{ij}$ is real and nonnegative  for $0 \leq h,i,j\leq D$.
The $q^h_{ij}$ are called the {\it Krein parameters}.
The graph $\Gamma$ is said to be {\it $Q$-polynomial}
(with respect to the given ordering $E_0, E_1, \ldots, E_D$
of the primitive idempotents)
whenever for $0 \leq h,i,j\leq D$, 
$q^h_{ij}= 0$
(resp. 
$q^h_{ij}\not= 0$) whenever one of $h,i,j$ is greater than
(resp. equal to) the sum of the other two \cite[p.~59]{bcn}.
See \cite{bannai,
caugh1,
caugh2,
aap1,
aap2} for more information on the $Q$-polynomial property.
From now on we assume that $\Gamma$ is 
$Q$-polynomial with respect to $E_0, E_1, \ldots, E_D$.

\medskip
\noindent
We  recall the dual Bose-Mesner algebra of $\Gamma.$
Fix a vertex $x \in X.$ We view $x$ as a ``base vertex.''
For 
$ 0 \le i \le D$ let $E_i^*=E_i^*(x)$ denote the diagonal
matrix in $\hbox{Mat}_X(\C)$ with $yy$ entry
\begin{equation}\label{DEFDEI}
{(E_i^*)_{yy} = \cases{1, & if $\partial(x,y)=i$\cr
0, & if $\partial(x,y) \ne i$\cr}} \qquad (y \in X).
\end{equation}
We call $E_i^*$ the  $i$th {\it dual idempotent} of $\Gamma$
 with respect to $x$ \cite[p.~378]{terwSub1}.
We observe that
(i) $\sum_{i=0}^D E_i^*=I$;
(ii) $\overline{E_i^*} = E_i^*$ $(0 \le i \le D)$;
(iii) $E_i^{*t} = E_i^*$ $(0 \le i \le D)$;
(iv) $E_i^*E_j^* = \delta_{ij}E_i^* $ $(0 \le i,j \le D)$.
By these facts 
$E_0^*,E_1^*, \ldots, E_D^*$ form a 
basis for a commutative subalgebra
$M^*=M^*(x)$ of 
$\hbox{Mat}_X(\C).$ 
We call 
$M^*$ the {\it dual Bose-Mesner algebra} of
$\Gamma$ with respect to $x$ \cite[p.~378]{terwSub1}.
For $0 \leq i \leq D$ let $A^*_i = A^*_i(x)$ denote the diagonal
matrix in 
 $\hbox{Mat}_X(\C)$
with $yy$ entry
$(A^*_i)_{yy}=\vert X \vert (E_i)_{xy}$ for $y \in X$.
Then $A^*_0, A^*_1, \ldots, A^*_D$ is a basis for $M^*$ 
\cite[p.~379]{terwSub1}.
Moreover
(i) $A^*_0 = I$;
(ii)
$\overline{A^*_i} = A^*_i \;(0 \le i \le D)$;
(iii) $A^{*t}_i = A^*_i  \;(0 \le i \le D)$;
(iv) $A^*_iA^*_j = \sum_{h=0}^D q_{ij}^h A^*_h \;( 0 \le i,j \le D)
$
\cite[p.~379]{terwSub1}.
We call 
 $A^*_0, A^*_1, \ldots, A^*_D$
the {\it dual distance matrices} of $\Gamma$ with respect to $x$.
We abbreviate
$A^*:=A^*_1$ 
and call this the {\it dual adjacency matrix} of $\Gamma$ with
respect to $x$.
The matrix $A^*$ generates $M^*$ \cite[Lemma 3.11]{terwSub1}.
By (\ref{SHF}) and since $A^*$ is real symmetric,
\begin{equation}
\label{SHF2}
\langle A^{*}u, v\rangle =\langle u, A^{*}v \rangle \qquad \qquad u,v \in V.
\end{equation}

\medskip
\noindent We recall the dual eigenvalues of $\Gamma$.
Since $E^*_0,E^*_1,\ldots,E^*_D$ form a basis for  
$M^*$ and since $A^*$ is real, there exist real scalars $\theta^*_0,\theta^*_1,
\ldots,\theta^*_D$ such that
$A^* = \sum_{i=0}^D \theta^*_iE^*_i$.
Observe that
$A^*E^*_i = E^*_iA^* =  \theta^*_iE^*_i$ for $0 \leq i \leq D$.
 We call $\theta^*_i$ the {\it dual eigenvalue}
of $\Gamma$ associated with $E^*_i$ $(0 \leq i\leq D)$. Observe that $\theta^*_0,\theta^*_1,\ldots,\theta^*_D$ are mutually
distinct since $A^*$ generates $M^*$.

\medskip
\noindent 
We recall the subconstituents of $\Gamma $.
From
(\ref{DEFDEI}) we find
\begin{equation}\label{DEIV}
E_i^*V = \mbox{span}\{\hat{y} \;|\; y \in X, \quad \partial(x,y)=i\}
\qquad (0 \le i \le D).
\end{equation}
By 
(\ref{DEIV})  and since
 $\{\hat{y}\;|\;y \in X\}$ is an orthonormal basis for $V$
 we find
\begin{eqnarray*}
\label{vsub}
V = E_0^*V+E_1^*V+ \cdots +E_D^*V \qquad \qquad 
{\rm (orthogonal\ direct\ sum}).
\end{eqnarray*}
For $0 \leq i \leq D$ the space $E^*_iV$ is the eigenspace
of $A^*$ associated with $\theta^*_i$.
We call $E_i^*V$ the {\it $i$th subconstituent} of $\Gamma$
with respect to $x$.

\medskip
\noindent
We recall the subconstituent algebra of $\Gamma $.
Let $T=T(x)$ denote the subalgebra of $\hbox{Mat}_X(\C)$ generated by 
$M$ and $M^*$. 
We call $T$ the {\it subconstituent algebra} 
(or {\it Terwilliger algebra}) of $\Gamma$ 
 with respect to $x$ \cite[Definition 3.3]{terwSub1}.
We observe that $T$ is generated by $A,A^*$.
We observe that $T$ has finite dimension. Moreover $T$ is 
semi-simple since it
is closed under the conjugate transponse map
\cite[p.~157]{CR}. 
See
\cite{curtin1,
curtin2,
go,
go2,
hobart,
tanabe,
terwSub1,
terwSub2,
terwSub3}
for more information on the subconstituent
algebra.

\medskip
\noindent 
For the rest of this paper we adopt the
following notational convention.

\begin{notation}
\label{setup}
\rm
We assume that $\Gamma=(X,R)$ is a  distance-regular graph with
diameter $D\geq 3$. We assume that $\Gamma$ is $Q$-polynomial
with respect to the ordering $E_0, E_1, \ldots, E_D$
of the primitive idempotents. We fix $x \in X$ 
and write $A^*=A^*(x)$, 
$E^*_i=E^*_i(x)$ $(0 \leq i \leq D)$,
$T=T(x)$. We abbreviate $V=\C^X$.
For notational convenience we define
$E_{-1}=0$, 
$E_{D+1}=0$ and 
$E^*_{-1}=0$, 
$E^*_{D+1}=0$.
\end{notation}

\noindent We finish this section with a comment.

\begin{lemma}
\label{lem:incl}
{\rm
\cite[Lemma 3.2]{terwSub1}}
With reference to Notation \ref{setup},
the following (i), (ii) hold for $0 \leq i \leq D$.
\begin{enumerate}
\item
$AE^*_iV \subseteq E^*_{i-1}V +E^*_iV+E^*_{i+1}V
$.
\item $A^*E_iV \subseteq E_{i-1}V +E_iV+E_{i+1}V
$.
\end{enumerate}
\end{lemma}

\section{The irreducible $T$-modules}

\noindent In this section 
we recall some useful results
on $T$-modules.

\medskip
\noindent
With reference to Notation \ref{setup},
by a {\it T-module}
we mean a subspace $W \subseteq V$ such that $BW \subseteq W$
for all $B \in T.$ 
 Let $W$ denote a $T$-module. Then $W$ is said
to be {\it irreducible} whenever $W$ is nonzero and $W$ contains 
no $T$-modules other than 0 and $W$.

\medskip
\noindent
Let $W$ denote a $T$-module and let 
$W'$ denote a  
$T$-module contained in $W$.
Then the orthogonal complement of $W'$ in $W$ is a $T$-module 
\cite[p.~802]{go2}.
It follows that each $T$-module
is an orthogonal direct sum of irreducible $T$-modules.
In particular $V$ is an orthogonal direct sum of irreducible $T$-modules.

\medskip
\noindent 
Let $W$ denote an irreducible $T$-module.
By the {\it endpoint} of $W$ we mean
$\mbox{min}\lbrace i |0\leq i \leq D, \; E^*_iW\not=0\rbrace $.
By the {\it diameter} of $W$ we mean
$ |\lbrace i | 0 \leq i \leq D,\; E^*_iW\not=0 \rbrace |-1 $.
By the {\it dual endpoint} of $W$ we mean
$\mbox{min}\lbrace i |0\leq i \leq D, \; E_iW\not=0\rbrace $.
By
the {\it dual diameter} of $W$ we mean
$ |\lbrace i | 0 \leq i \leq D,\; E_iW\not=0 \rbrace |-1 $.
The diameter of $W$ is  equal to the dual diameter of
$W$
\cite[Corollary 3.3]{aap1}.

\begin{lemma}
{\rm
\cite[Lemma 3.4, Lemma 3.9, Lemma 3.12]{terwSub1}}
\label{lem:basic}
With reference to Notation \ref{setup},
let $W$ denote an irreducible $T$-module with endpoint $\rho$,
dual endpoint $\tau$, and diameter $d$.
Then $\rho,\tau,d$ are nonnegative integers such that $\rho+d\leq D$ and
$\tau+d\leq D$. Moreover the following (i)--(iv) hold.
\begin{enumerate}
\item 
$E^*_iW \not=0$ if and only if $\rho \leq i \leq \rho+d$ 
$ \quad (0 \leq i \leq D)$.
\item
$W = \sum_{h=0}^{d} E^*_{\rho+h}W$ \qquad \rm({orthogonal direct sum}).
\item 
$E_iW \not=0$ if and only if $\tau \leq i \leq \tau+d$
$ \quad (0 \leq i \leq D)$.
\item
$W = \sum_{h=0}^{d} E_{\tau+h}W$ \qquad \rm({orthogonal direct sum}).
\end{enumerate}
\end{lemma}

\begin{lemma}
{\rm
\cite[Lemma 3.2]{ds}}
\label{lem:wincl}
With reference to Notation \ref{setup},
let $W$ denote an irreducible $T$-module with
endpoint $\rho$,
dual endpoint $\tau$, and diameter $d$.
Then the following (i), (ii) hold for $0 \leq i \leq d$.
\begin{enumerate}
\item
$AE^*_{\rho+i}W \subseteq E^*_{\rho+i-1}W +E^*_{\rho+i}W+E^*_{\rho+i+1}W
$.
\item $A^*E_{\tau+i}W \subseteq E_{\tau+i-1}W +E_{\tau+i}W+E_{\tau+i+1}W
$.
\end{enumerate}
\end{lemma}

\begin{remark} 
\label{rem:tdpair}
\rm
With reference to Notation \ref{setup},
let $W$ denote an irreducible $T$-module.
Then $A$ and  $A^*$ act on $W$ as a tridiagonal pair
in the sense of Definition \ref{def:nonthinlp}.
This follows from
Lemma \ref{lem:basic},
Lemma 
\ref{lem:wincl}, and since $A,A^*$ together generate
$T$. 
\end{remark}

\begin{lemma}
\label{lem:wsplit}
With reference to Notation \ref{setup}, let $W$ 
denote an irreducible $T$-module with endpoint $\rho$,
dual endpoint $\tau$, and diameter $d$.
Then for $\mu, \nu \in \lbrace \downarrow,\uparrow \rbrace$ we have
\begin{eqnarray}
\label{eq:splitd}
W = \sum_{h=0}^d W^{\mu \nu}_h \qquad (\rm{direct\ sum}),
\end{eqnarray}
where for $0 \leq h\leq d$,
\begin{eqnarray*}
W^{\downarrow \downarrow}_h &=& (E^*_\rho W+\cdots+ E^*_{\rho+h}W)\cap (E_\tau W+\cdots + E_{\tau+d-h}W),
\label{eq:whdef1}
\\
W^{\uparrow \downarrow}_h &=& (E^*_{\rho+d-h} W+\cdots+ E^*_{\rho+d}W)\cap (E_\tau W+\cdots + E_{\tau+d-h}W),
\label{eq:whdef2}
\\
W^{\downarrow \uparrow}_h &=& (E^*_\rho W+\cdots+ E^*_{\rho+h}W)\cap (E_{\tau+h}W+\cdots + E_{\tau+d}W),
\label{eq:whdef3}
\\
W^{\uparrow \uparrow}_h &=& (E^*_{\rho+d-h} W+\cdots+ E^*_{\rho+d}W)\cap (E_{\tau+h}W+\cdots + E_{\tau+d}W).
\label{eq:whdef4}
\end{eqnarray*}
\end{lemma}

\noindent 
{\it Proof:} Immediate from Lemma \ref{lem:spl} and Remark \ref{rem:tdpair}.
\hfill $\Box$ \\

\noindent We remark that the sum (\ref{eq:splitd}) is not orthogonal in general. However we do have the following result.

\begin{lemma}
\label{lem:wsplitorg}
With reference to Notation \ref{setup}, let $W$ 
denote an irreducible $T$-module with diameter $d$.
Then the following (i), (ii) hold for $0 \leq h,\ell \leq d$ such that $h+\ell \neq d$.
\begin{enumerate}
\item The subspaces $W^{\downarrow \downarrow}_h$ and $W^{\uparrow \uparrow}_{\ell}$ are orthogonal with respect to the standard Hermitian form.
\item The subspaces $W^{\downarrow \uparrow}_h$ and $W^{\uparrow \downarrow}_{\ell}$ are orthogonal with respect to the standard Hermitian form.
\end{enumerate}
\end{lemma}

\noindent 
{\it Proof:} Combine Theorem \ref{thm:ortho2}, (\ref{SHF1}), (\ref{SHF2}), Remark \ref{rem:tdpair}, and Lemma \ref{lem:wsplit}.
\hfill $\Box$

\section{The split decompositions of the standard module}

\noindent In this section we recall the four split decompositions for the standard module and discuss their basic properties.

\begin{definition}
{\rm
\cite[Definition 10.1]{drgqtet}}
\label{def:updown}
\rm
With reference to Notation \ref{setup}, for $-1\leq i,j\leq D$
we define
\begin{eqnarray*}
V^{\downarrow \downarrow}_{i,j} &=& 
(E^*_0V+\cdots+E^*_iV)\cap (E_0V+\cdots+E_jV),
\\
V^{\uparrow \downarrow}_{i,j} &=& 
(E^*_DV+\cdots+E^*_{D-i}V)\cap (E_0V+\cdots+E_jV),
\\
V^{\downarrow \uparrow}_{i,j} &=& 
(E^*_0V+\cdots+E^*_iV)\cap (E_DV+\cdots+E_{D-j}V),
\\
V^{\uparrow \uparrow}_{i,j} &=& 
(E^*_DV+\cdots+E^*_{D-i}V)\cap (E_DV+\cdots+E_{D-j}V).
\end{eqnarray*}
In each of the above four equations we interpret the
right-hand side to be 0 if $i=-1$ or $j=-1$. 
\end{definition}

\begin{definition}
{\rm
\cite[Definition 10.2]{drgqtet}}
\label{def:vtilde}
\rm
With reference to Notation \ref{setup} and Definition \ref{def:updown}, for $\mu, \nu \in \lbrace \downarrow,
\uparrow \rbrace$ and  $0 \leq i,j\leq D$ we have
$
V^{\mu \nu}_{i-1,j} \subseteq V^{\mu \nu}_{i,j}$
and
$
V^{\mu \nu}_{i,j-1}  \subseteq V^{\mu \nu}_{i,j}
$. Therefore
\begin{eqnarray*}
V^{\mu \nu}_{i-1,j}+
V^{\mu \nu}_{i,j-1} \subseteq V^{\mu \nu}_{i,j}.
\end{eqnarray*}
Referring to the above inclusion, we define ${\tilde V}^{\mu \nu}_{i,j}$
 to be the orthogonal complement of the left-hand side in the
 right-hand side; that
is
\begin{eqnarray*}
 {\tilde V}^{\mu \nu}_{i,j}=(
V^{\mu \nu}_{i-1,j}+
V^{\mu \nu}_{i,j-1})^\perp \cap V^{\mu \nu}_{i,j}.
\end{eqnarray*}
\end{definition}

\begin{lemma}
\label{lem:drsum}
{\rm
\cite [Lemma 10.3]{drgqtet}}
With reference to Notation \ref{setup} and Definition \ref{def:vtilde}, the following holds 
 for $\mu, \nu \in \lbrace \downarrow, \uparrow \rbrace$:
\begin{eqnarray}
V = \sum_{i=0}^D \sum_{j=0}^D {\tilde V}^{\mu \nu}_{i,j}
\qquad \qquad \rm({direct\ sum}).
\label{eq:vsplt}
\end{eqnarray}
\end{lemma}

\begin{definition}
\label{def:splitdv}
\rm
We call the sum 
(\ref{eq:vsplt})
the $(\mu,\nu)$--{\it split decomposition} of $V$ with respect to $x$.
\end{definition}

\begin{remark}
\rm
The decomposition (\ref{eq:vsplt}) is
not orthogonal in general.
\end{remark}

\begin{lemma}
\label{lem:DU}
With reference to Notation \ref{setup}, let $W$ denote an irreducible $T$-module with endpoint $\rho$, dual endpoint $\tau$, and diameter $d$. Then for $0 \leq h\leq d$ and $0 \leq i,j\leq D$ the following (i)--(iv) hold.
\begin{enumerate}
\item $W^{\downarrow \downarrow}_h\subseteq {{\tilde V}^{\downarrow \downarrow}_{i,j}}$ if and only if $i=\rho +h$ and $j=\tau+d-h$.
\item $W^{\uparrow \downarrow}_h\subseteq {{\tilde V}^{\uparrow \downarrow}_{i,j}}$ if and only if $i=D-\rho-d+h$ and $j=\tau +d-h$.
\item $W^{\downarrow \uparrow}_h\subseteq {{\tilde V}^{\downarrow \uparrow}_{i,j}}$ if and only if $i=\rho +h$ and $j=D-\tau-h$.
\item $W^{\uparrow \uparrow}_h\subseteq {{\tilde V}^{\uparrow \uparrow}_{i,j}}$ if and only if $i=D-\rho-d+h$ and $j=D-\tau-h$.
\end{enumerate}
\end{lemma}

\noindent {\it Proof:}
Immediate from \cite[Lemma 11.4]{drgqtet} and (\ref{eq:vsplt}).
\hfill $\Box$ \\

\begin{lemma}
\label{lem:DUsum}
With reference to Notation \ref{setup}, fix an orthogonal direct sum decomposition of the standard module $V$ of $\Gamma$ into irreducible $T$-modules:
\begin{eqnarray}
V = \sum_{W} {W}.
\label{eq:orthods}
\end{eqnarray}
Then the following (i)--(iv) hold for $0 \leq i,j\leq D$.
\begin{enumerate}
\item ${{\tilde V}^{\downarrow \downarrow}_{i,j}} =\sum{W^{\downarrow \downarrow}_h}$, where the sum is over all ordered pairs $(W,h)$ such that $W$ is assumed in (\ref {eq:orthods}) with endpoint $\rho \leq i$, dual endpoint $\tau=i+j-\rho-d$, diameter $d \geq i-\rho$, and $h=i- \rho$.

\item ${{\tilde V}^{\uparrow \downarrow}_{i,j}} =\sum{{W^{\uparrow \downarrow}_h}}$, where the sum is over all ordered pairs $(W,h)$ such that $W$ is assumed in (\ref {eq:orthods}) with endpoint $\rho \leq D-i$, dual endpoint $\tau=i+j+\rho-D$, diameter $d \geq D-\rho-i$, and $h=\rho+ d-D+i$.

\item ${{\tilde V}^{\downarrow \uparrow}_{i,j}} =\sum{{W^{\downarrow \uparrow}_h}}$, where the sum is over all ordered pairs $(W,h)$ such that $W$ is assumed in (\ref {eq:orthods}) with endpoint $\rho \leq i$, dual endpoint $\tau=\rho+D-i-j$, diameter $d \geq i-\rho$, and $h=i-\rho$.

\item ${{\tilde V}^{\uparrow \uparrow}_{i,j}} =\sum{{W^{\uparrow \uparrow}_h}}$, where the sum is over all ordered pairs $(W,h)$ such that $W$ is assumed in (\ref {eq:orthods}) with endpoint $\rho \leq D-i$, dual endpoint $\tau=2D-\rho-d-i-j$, diameter $d \geq D-\rho-i$, and $h=\rho+ d-D+i$.
\end{enumerate}
\end{lemma}

\noindent {\it Proof:}
(i) For $0 \leq i,j\leq D$ define
\begin{eqnarray} 
v_{i,j} &=& \sum{W^{\downarrow \downarrow}_h},
\label{eq:ddsum}
\end{eqnarray}
where the sum is over all ordered pairs $(W,h)$ such that $W$ is assumed in (\ref {eq:orthods}) with endpoint $\rho \leq i$, dual endpoint $\tau=i+j-\rho-d$, diameter $d \geq i-\rho$, and $h=i- \rho$. We show that ${{\tilde V}^{\downarrow \downarrow}_{i,j}}=v_{i,j}$. We first show that ${{\tilde V}^{\downarrow \downarrow}_{i,j}} \supseteq v_{i,j}$. Let $W^{\downarrow \downarrow}_h$ denote one of the terms in the sum on the right in (\ref {eq:ddsum}). We show that $W^{\downarrow \downarrow}_h$ is contained in ${{\tilde V}^{\downarrow \downarrow}_{i,j}}$. Let $\rho, \tau,d$ denote the endpoint, dual endpoint, and diameter of $W$, respectively. By construction $\tau=i+j-\rho-d$ and $h=i- \rho$. Subtracting the second equation from the first equation we find $j=\tau+d-h$. Now $W^{\downarrow \downarrow}_h$ is contained in ${{\tilde V}^{\downarrow \downarrow}_{i,j}}$ by Lemma \ref{lem:DU}(i). We have now shown that ${{\tilde V}^{\downarrow \downarrow}_{i,j}} \supseteq v_{i,j}$. We can now easily show that ${{\tilde V}^{\downarrow \downarrow}_{i,j}}=v_{i,j}$. Expanding the sum (\ref {eq:orthods}) using Lemma \ref {lem:wsplit} we get
\begin{eqnarray*} 
V &=& \sum_{W} {W} \qquad \qquad (\mbox{direct sum})\\
  &=& \sum_{W} \sum_{h} W^{\downarrow \downarrow}_h \qquad \qquad (\mbox{direct sum}),
\end{eqnarray*}
where the second sum is over the integer $h$ from 0 to the diameter of $W$. In the above sum we change the order of summation to get
\begin{eqnarray*} 
V &=& \sum_{i=0}^D \sum_{j=0}^D \sum{W^{\downarrow \downarrow}_h} \qquad \qquad (\mbox{direct sum}),
\end{eqnarray*}
where the third sum is over all ordered pairs $(W,h)$ such that $W$ is assumed in (\ref {eq:orthods}) with endpoint $\rho \leq i$, dual endpoint $\tau=i+j-\rho-d$, diameter $d \geq i-\rho$, and $h=i- \rho$. In other words,
\begin{eqnarray*} 
V &=& \sum_{i=0}^D \sum_{j=0}^D {v}_{i,j} \qquad \qquad (\mbox{direct sum}).
\end{eqnarray*}
By this, (\ref{eq:vsplt}), and since ${{\tilde V}^{\downarrow \downarrow}_{i,j}} \supseteq v_{i,j}$ for $0 \leq i,j\leq D$, we find ${{\tilde V}^{\downarrow \downarrow}_{i,j}}=v_{i,j}$ for $0 \leq i,j\leq D$.\\
(ii), (iii), (iv) Similar to the proof of (i).
\hfill $\Box$ \\

\noindent 
Now we have the main result.

\begin{theorem}
\label{thm:duality}
With reference to Notation \ref{setup} and Definition \ref{def:vtilde}, the following (i), (ii) hold for $0 \leq i,j,r,s\leq D$.
\begin{enumerate}
\item ${{\tilde V}^{\downarrow \downarrow}_{i,j}}$ and ${{\tilde V}^{\uparrow \uparrow}_{r,s}}$ are orthogonal unless $i+r=D$ and $j+s=D$.
\item ${{\tilde V}^{\downarrow \uparrow}_{i,j}}$ and ${{\tilde V}^{\uparrow \downarrow}_{r,s}}$ are orthogonal unless $i+r=D$ and $j+s=D$.
\end{enumerate}
\end{theorem}

\noindent {\it Proof:}
(i) Assume that $i+r \neq D$ or $j+s \neq D$. We show that ${{\tilde V}^{\downarrow \downarrow}_{i,j}}$ and ${{\tilde V}^{\uparrow \uparrow}_{r,s}}$ are orthogonal. To do this we will use Lemma \ref{lem:DUsum}(i),(iv). 
Let $W^{\downarrow \downarrow}_h$ (resp. ${W'}^{\uparrow \uparrow}_{h'}$) denote one of the terms in the sum in Lemma \ref{lem:DUsum}(i) (resp. Lemma \ref{lem:DUsum}(iv)). We show that $W^{\downarrow \downarrow}_h$ and ${W'}^{\uparrow \uparrow}_{h'}$ are orthogonal. There are two cases to consider. First assume that $W \neq W'$. Then $W$ and $W'$ are orthogonal so $W^{\downarrow \downarrow}_h$ and ${W'}^{\uparrow \uparrow}_{h'}$ are orthogonal. Next assume that $W=W'$. Let $\rho, \tau,d$ denote the corresponding endpoint, dual endpoint, and diameter. By Lemma \ref{lem:DUsum}(i),
\begin{equation}
\label{main1}
\tau=i+j-\rho-d,\qquad \qquad h=i- \rho. 
\end{equation}
By Lemma \ref{lem:DUsum}(iv),\\
\begin{equation}
\label{main2}
\tau=2D-\rho-d-r-s, \qquad \qquad h'=\rho+ d-D+r. 
\end{equation}
Adding the equations on the right in (\ref {main1}), (\ref {main2}) we get
\begin{equation}
\label{main3}
i+r-D=h+h'-d. 
\end{equation}
Subtracting the equation on the left in (\ref {main1}) from the equation on the left in (\ref {main2}) and evaluating the result using (\ref {main3}) we get
\begin{equation}
\label{main4}
j+s-D=d-h-h'. 
\end{equation}
By (\ref {main3}), (\ref {main4}) and since $i+r \neq D$ or $j+s \neq D$ we find $h+h' \neq d$.
Now $W^{\downarrow \downarrow}_h$ and ${W'}^{\uparrow \uparrow}_{h'}$ are orthogonal by Lemma \ref{lem:wsplitorg}(i).\\
(ii) Similar to the proof of (i).
\hfill $\Box$ \\

\begin{corollary}
With reference to Notation \ref{setup} and Definition \ref{def:vtilde}, the following (i), (ii) hold for $0 \leq i,j\leq D$.
\begin{enumerate}
\item {\rm{dim}}$\,{{\tilde V}^{\downarrow \downarrow}_{i,j}}$ $=$ {\rm{dim}}$\,{{\tilde V}^{\uparrow \uparrow}_{D-i,D-j}}$.
\item {\rm{dim}}$\,{{\tilde V}^{\downarrow \uparrow}_{i,j}}$ $=$ {\rm{dim}}$\,{{\tilde V}^{\uparrow \downarrow}_{D-i,D-j}}$.
\end{enumerate}
\end{corollary}

\noindent {\it Proof:}
Immediate from Theorem \ref{thm:duality} and elementary linear algebra.
\hfill $\Box$ \\

\bigskip
\noindent{\Large\bf Acknowledgements}

\medskip
\noindent 
The author would like to thank Professor Paul M. Terwilliger for his valuable ideas and comments.

\noindent Joohyung Kim \\
 Department of Mathematics \\ 
University of Wisconsin \\
480 Lincoln Drive \\
Madison Wisconsin \\
53706-1388 USA\\
Email: jkim@math.wisc.edu\\

\end{document}